\documentclass[11pt]{article}

\usepackage{geometry}                
\usepackage[parfill]{parskip}  
\usepackage{graphicx}
\usepackage{amssymb}
\usepackage{epstopdf}
\usepackage{amsmath,amsthm}
\usepackage{tabularx}
\usepackage{algorithmicx}
\usepackage[ruled]{algorithm}
\usepackage{algpseudocode}
\usepackage{url}
\usepackage{authblk}
\usepackage{todonotes}
\usepackage{color,soul}

\providecommand{\keywords}[1]{\textbf{{Keywords }} #1}
\providecommand{\subjclass}[1]{\textbf{{Mathematics subject classification }} #1}
\algrenewcommand\algorithmicindent{1.0em}

\newtheorem{theorem}{Theorem}

\title{High-Dimensional Sparse Fourier Algorithms}

\author{%
Bosu Choi\thanks{The Oden Institute for Computational Engineering and Sciences, University of Texas at Austin \texttt{choibosu@utexas.edu}}
\and
Andrew Christlieb  \thanks{Computational Mathematics Science and Engineering, Michigan State University \texttt{christli@msu.edu}}%
\and
Yang Wang\thanks{Department of Mathematics, The Hong Kong University of Science and Technology  \texttt{yangwang@ust.hk‎} }
}

\date{}

\begin{document}

\maketitle

\begin{abstract}
In this paper, we discuss the development of a sublinear sparse Fourier algorithm for high-dimensional data. In \lq\lq Adaptive Sublinear Time Fourier Algorithm" by D. Lawlor, Y. Wang and  A. Christlieb (2013) \cite{lawlor2013adaptive}, an efficient algorithm with $\Theta(k\log k)$ average-case runtime and $\Theta(k)$ average-case sampling complexity for the one-dimensional sparse FFT was developed for signals of bandwidth $N$, where $k$ is the  number of significant modes such that $k\ll N$.

In this work we develop an efficient algorithm for sparse FFT for higher dimensional signals, extending some of the ideas in  \cite{lawlor2013adaptive}. Note a higher dimensional signal can always be unwrapped into a one dimensional signal, but when the dimension gets large, unwrapping a higher dimensional signal into a one dimensional array is far too expensive to be realistic. Our approach here introduces two new concepts: \lq\lq partial unwrapping'' and \lq\lq tilting''. These two ideas allow us to efficiently compute the sparse FFT of higher dimensional signals.
\end{abstract}



\section{INTRODUCTION} \label{INTRODUCTION}
\setcounter{equation}{0}

  As the size and dimensionality of data sets in science and engineering  grow larger and larger, it is necessary to develop efficient tools to analyze them \cite{greene2015understanding, LSST}. One of the best known and most frequently-used tools is the Fast Fourier Transform (FFT). However, in the case that the bandwidth $N$ of frequencies is large, the sampling size becomes large, as dictated by the Shannon-Nyquist sampling theorem. Specifically, the runtime complexity is $\mathcal{O}(N{\log}N)$ and the number of samples is $\mathcal{O}(N)$. This issue is only exacerbated in the $d$-dimensional setting, where the runtime complexity is $\mathcal{O}(N^d{\log}N^d )$ and the number of samples is $\mathcal{O}(N^d)$ if we assume the dimension is $d$ and the bandwidth in each dimension is $N$.  Due to this \lq\lq curse of dimensionality'', many higher dimensional problems of interest are beyond current computational capabilities of the traditional  FFT. Moreover, in the sparse setting where the number of significant frequencies $k$ is small, it is computationally wasteful to compute all $N^d$ coefficients.  In such a setting we refer to the problem as being \lq\lq sparse''. For sparse problems, the idea of sublinear sparse Fourier transforms was introduced  \cite{gilbert2002near, gilbert2005improved, hassanieh2012nearly, hassanieh2012simple, iwen2010combinatorial, lawlor2013adaptive, christlieb2016multiscale, plonka2016, plonka2018, iwen2013improved, 10.3389/fams.2016.00001}. These methods greatly reduce the runtime and sampling complexity of the FFT in the sparse setting. The methods were primarily  designed for the one dimensional setting.

  The  first sparse Fourier algorithm  was proposed in \cite{gilbert2002near}. It introduced a randomized algorithm with  $\mathcal{O}(k^2{\log}^c N)$ runtime and $\mathcal{O}(k^2{\log}^c N)$ samples where $c$ is a positive number that varies depending on the trade-off between efficiency and accuracy. An algorithm with improved runtime $\mathcal{O}(k{\log}^c N)$ and samples $\mathcal{O}(k{\log}^c N)$ was given in \cite{gilbert2005improved}. The algorithms given in \cite{hassanieh2012nearly} and \cite{hassanieh2012simple} achieved $\mathcal{O}(k{\log} N{\log} {N}/{k})$ average-case runtime and gave empirical results. The algorithms in \cite{gilbert2002near, gilbert2005improved, hassanieh2012nearly, hassanieh2012simple} are all randomized. The first deterministic algorithm using a combinatorial approach was introduced in \cite{iwen2010combinatorial}. In \cite{lawlor2013adaptive}, another deterministic algorithm was given whose procedure recognizes frequencies in a similar manner to \cite{hassanieh2012nearly}.  The two methods in \cite{lawlor2013adaptive, hassanieh2012nearly} were published at the same time and both use the idea of working with two sets of samples, one at $\mathcal{O}(k)$  points and the second at the same $\mathcal{O}(k)$ points plus a small shift. The ratio of the FFT of the two sets of points, plus extra machinery, lead to fast deterministic algorithms.  The first deterministic algorithm \cite{iwen2010combinatorial} has $\mathcal{O}(k^2{\log}^4 N)$ runtime and sampling complexity, and the second one \cite{lawlor2013adaptive} has $\mathcal{O}(k{\log}k)$ average-case runtime and $\mathcal{O}(k)$ sampling complexity. Later, \cite{christlieb2016multiscale} introduced modified methods for noisy data with $\mathcal{O}(k{\log}k{\log}N/k)$ average-case runtime and $\mathcal{O}(k{\log}N/k)$ sampling complexity. Also, there is a method under the assumption of block-structured sparsity in \cite{bittens2017deterministic}. Our method, discussed throughout this paper builds on the method presented in \cite{lawlor2013adaptive}.
  Moreover, extending the one-dimensional fully discrete sparse Fourier transforms from  \cite{merhi2017new} into a high-dimensional setting by using our method would be interesting future work. 

 The methods introduced in the previous paragraph are for one-dimensional data. In \cite{ghazi2013sample}, practical algorithms for data in two dimensions were given for the first time. In \cite{KAMMERER2015543}, high dimensional sparse FFT was introduced using the rank-1 lattice sampling, which shows the numerical results up to dimension 10. In this paper,  we develop algorithms designed for higher dimensional data, which is effective even for dimensions in the hundreds and thousands. To achieve our goal, our approach must address the worst case scenario presented in \cite{ghazi2013sample}. We can find a variety of data sets in multiple dimensions that we want to analyze. A relatively low-dimensional example is MRI data, which is three dimensional. However, when we designed the method in this paper, we had much higher dimensional problems in mind, such as some astrophysical data, e.g., the  Sloan Digital Sky Survey and Large Synoptic Survey Telescope \cite{SDSS, LSST}. They produce tera- or peta-bytes of imaging and spectroscopic data in very high dimensions.  Due to the computational effort of a multi-dimensional FFT, spectral analysis of  this high dimensional data  necessitates a multidimensional sparse fast Fourier transform. Further, given the massive size of data sets in some current and future problems in science and engineering, it is anticipated that the development of such an efficient algorithm will play an important role in the analysis of these types of data.

  It is not straightforward to extend one dimensional sparse Fourier transform algorithms to multiple dimensions. We face several obstacles. First,  we do not have an efficient FFT for multidimensional problems much higher than three. Using projections onto lower-dimensional spaces solves this problem. However, like all projection methods for sparse FFT, one needs to match frequencies from one projection with those from another projection. This {\em registration} problem is one of the big challenges in the one dimensional sparse FFT. An equally difficult challenge is that different frequencies may be projected into the same frequency ({\em the collision problem}). All projection methods for sparse FFT primarily aim to overcome these two challenges. In higher dimensional sparse FFT, these problems become even more challenging as now we are dealing with frequency vectors, not just scalar frequencies.

  As a first step to our goal of a high dimensional sparse FFT, this paper addresses the case for continuous data without noise in a high dimensional setting. We introduce effective methods to address the registration and the collision problems. In particular, we introduce a novel {\em partial unwrapping} technique that is shown to be highly effective in reducing the registration and collision complexity while maintains the sublinear runtime efficiency.
  We shall show that we can achieve $\Theta(dk{\log}k)$ average-case computational complexity and $\Theta(dk)$ average-case sampling complexity. In Section \ref{5}, we present as examples computational results for sparse FFT where the dimensions are 100 and 1000 respectively. For comparison, the traditional $d$-dimensional FFT requires $\mathcal{O}(N^d{\log}N^d)$ time complexity and $\mathcal{O}(N^d)$ sampling complexity, which is impossible to implement on any computers today.

\section{PRELIMINARIES} \label{PRELIMINARIES}
\setcounter{equation}{0}
\subsection{Review of the One-Dimensional Sublinear Sparse Fourier Algorithms} \label{2.1}
\begin{figure}[ht]
\centering
\includegraphics[width=0.9\textwidth, angle=0]{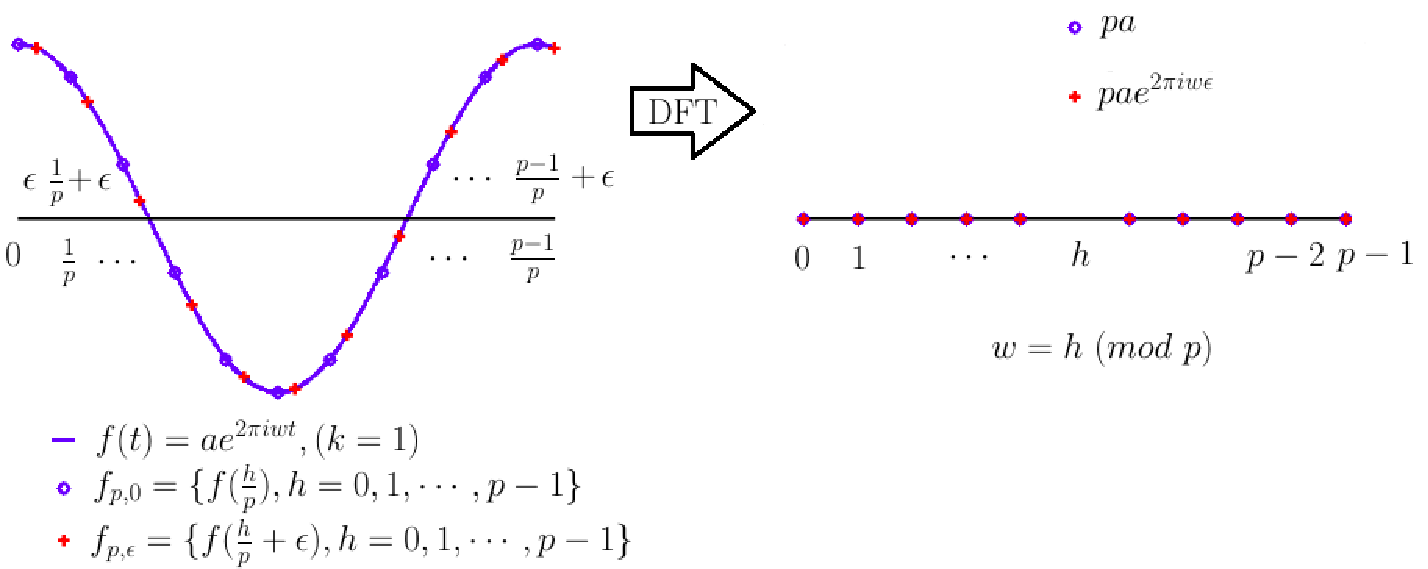}
\caption{Process of 1D sublinear sparse Fourier algorithm}\label{1D}
\label{fig1}
\end{figure}

 The one-dimensional sublinear sparse Fourier algorithm inspiring our method was developed in \cite{lawlor2013adaptive}. We briefly introduce the idea and notation of the algorithm before developing the multidimensional ones throughout this paper. We assume a function $f:[0,1)\rightarrow {\mathbb C}$ with sparsity $k$ as the following,
\begin{equation}\label{(1)}
	f(t) ~=~ \sum_{j=1}^{k}a_je^{2\pi i w_j t}
\end{equation}
with bandwidth $N$, i.e., frequency $w_j$ belongs to $[-N/2, N/2)\cap {\mathbb Z}$ and corresponding nonzero coefficient $a_j$ is in ${\mathbb C}$ for all $j$. We can consider it as a periodic function over ${\mathbb R}$ instead of $[0,1)$. The goal of the algorithm is to recover all coefficients $a_j$ and frequencies $w_j$ so that we can reconstruct the function $f$. This algorithm is called the \lq\lq phase-shift" method since it uses equi-spaced samples from the function and those at positions shifted by a small positive number $\epsilon$. To verify that the algorithm correctly finds the frequencies in the bandwidth $N$, $\epsilon$ should be strictly no bigger than $1/N$. We denote a sequence of samples shifted by $\epsilon$ with sampling rate $1/p$, where $p$ is a prime number, as
\begin{equation}\label{(2)}
	{\mathbf f}_{p, \epsilon} ~=~ \left( f(0 + \epsilon), f\left(\frac{1}{p} + \epsilon\right), f\left(\frac{2}{p} +\epsilon\right), f\left(\frac{3}{p}+ \epsilon\right), \cdots, f\left(\frac{p-1}{p} + \epsilon\right)\right).
\end{equation}
We skip much of the details here. In a nutshell, by choosing $p$ slightly larger than $k$ is enough to make the algorithm work. In \cite{lawlor2013adaptive} $p$ is set to be roughly $5k$, which is much smaller than the Nyquist rate $N$. Discrete Fourier transform (DFT) $\mathcal{F}$ of a sequence ${\mathbf b}=\left( b[0], b[1], \cdots, b[M-1] \right)$ of $M\in \mathbb{N}$ elements is defined as 
\begin{equation}
\mathcal{F}({\mathbf b})[h]:=\sum_{j=0}^{M-1} b[j] e^{-2 \pi i h  \cdot \frac{j}{M}}
\end{equation}
for $h=0,1,\cdots, M-1$.
If DFT is applied to the sample sequence ${\mathbf f}_{p, \epsilon}$, then the $h$-th element of its result is the following
\begin{equation}\label{(3)}
	{\mathcal F}({\mathbf f}_{p,\epsilon})[h] ~=~  p \sum_{w_j = h (\bmod p)} a_j e^{2 \pi i \epsilon w_j}
\end{equation}
where $h = 0, 1, \dots, p-1$. If there is only one frequency $w_j$ congruent to $h$ modulo $p$,
\begin{equation}\label{(4)}
	{\mathcal F}({\mathbf f}_{p,\epsilon})[h] ~=~ pa_je^{2 \pi i \epsilon w_j}.
\end{equation}
By putting $0$ instead of $\epsilon$, we can get unshifted samples ${\mathbf f}_{p,0}$ and applying the DFT gives
\begin{equation}\label{(5)}
	{\mathcal F}({\mathbf f}_{p,0})[h] ~=~ pa_j.
\end{equation}
This process so far is visualized in the Figure \ref{1D}. As long as there is no collision of frequencies with modulo $p$, we can find frequencies and their corresponding coefficients by the following computation
\begin{align}
	w_j &~=~ \frac{1}{2\pi \epsilon}{\rm Arg}\Big( \frac{{\mathcal F}({\mathbf f}_{p,\epsilon})[h]}      {{\mathcal F}({\mathbf f}_{p,0})[h]} \Big) ,\nonumber\\
	a_j &~=~ \frac{1}{p}{\mathcal F}({\mathbf f}_{p,0})[h],\label{(7)}
\end{align}
where the function $\lq\lq {\rm Arg}"$ gives us the argument falling into $[-\pi, \pi).$ Note that $w_j$ should be the only frequency congruent to $h$ modulo $p$, i.e., $w_j$ has no collision with other frequencies modulo $p$. The test to determine whether collision occurs or not is
\begin{equation}\label{(8)}
	\frac{\vert{\mathcal F}({\mathbf f}_{p,\epsilon})[h]\vert}{\vert{\mathcal F}({\mathbf f}_{p,0})[h]\vert}~=~ 1.
\end{equation}
The equality above holds when there is no collision. If there is a collision, the equality does not hold for almost all $\epsilon$ since the test fails to predict a collision for finite number of $\epsilon$ \cite{lawlor2013adaptive}. Further, it is also shown in the same paper that for any $\epsilon=\frac{a}{b}$ with $a, b$ coprime and $b\geq 2N$, the equality in (\ref{(8)}) does not hold at least once in a particular finite number of choice of $\epsilon$ unless there is no collision. In practical implementations, we choose $\epsilon$ to be ${1}/{cN}$ for some positive integer $c \geq 1$ and allow some small difference $\tau$ between the left and right sides of  (\ref{(8)}) where $\tau$ is very small positive number. For the choice of $\tau$, $p/N$ works well in our numerical experiments.
 
  The above process is one loop of the algorithm with a prime number $p$. To explain it from a different view, we can imagine that there are $p$ bins. Then we sort all frequencies into these bins according to their remainder modulo $p$. If there are more than one frequencies in one bin, then a collision happens. If there is only one frequency, then there is no collision. To determine whether a collision occurs, we use the above test. In the case where the test fails, i.e., the ratio is not $1$, we need to use another prime number $p'$. Thus we re-sort the frequencies into $p'$ bins by their remainder modulo $p'$. Even if two frequencies collide modulo $p$, it is likely that they do not collide modulo $p'$. Particularly, the Chinese Remainder Theorem guarantees that with a finite set of prime numbers, $\{p_{\ell}\}$, any frequency within the bandwidth $N$ can be uniquely identified, given $\prod_{\ell} p_{\ell} \geq N$. Algorithmically, for each loop, we choose a different prime number $p'$ and repeat equations (\ref{(2)})-(\ref{(8)}) with $p$ replaced by $p'$. In this way we can recover all $a_j$ and $w_j$ in  runtime $\Theta(k{\log} k)$ using $\Theta(k)$ samples on average. The overall code is shown in Algorithm \ref{Algorithm1} referred from \cite{lawlor2013adaptive}.

\vspace{-0.2cm}
\alglanguage{pseudocode}
\begin{algorithm}[ht]
\small
\caption{Phaseshift}
\label{Algorithm1}
\begin{algorithmic}[1]
\Procedure{$\mathbf{Phaseshift}$}{}\\
{\textbf{Input:}}{$f, c, k, N, \epsilon$}\\
{\textbf{Output:}}{$R$}
\State $R \gets \emptyset$
\State $i \gets 1$
\While { $|R|<k$}
    \State $ k^{\ast} \gets k - |R|$
    \State $p \gets {\it{i}\text{-th}}$ prime number $\geq c k^{\ast}$
    \State $g({{t}})=\sum_{({w},a_{w})\in R} a_{w} e^{2 \pi i {w} {t}}$

           \For {$h = 1 \to p$}
              \State ${f}_{p,\epsilon}[h] = f(\frac{h}{p}) - g(\frac{h}{p})  $
             \State ${f}_{p,0}[h] = f(\frac{h}{p}+\epsilon) - g(\frac{h}{p}+\epsilon)$
           \EndFor
           \State $\mathcal{F}({\mathbf f}_{p, \epsilon})=FFT({\mathbf f}_{p, \epsilon})$
           \State $\mathcal{F}({\mathbf f}_{p, 0})=FFT({\mathbf f}_{p, 0})$
           \State $\mathcal{F}^{sort}({\mathbf f}_{p, 0})=SORT(\mathcal{F}({\mathbf f}_{p, 0}))$

        \For {$h = 1 \to k^{\ast}$}
              \If {$\Big|\frac{|\mathcal{F}^{sort}({\mathbf f}_{p, 0})[h]|}{|\mathcal{F}^{sort}({\mathbf f}_{p, \epsilon})[h]|}-1\Big|<\epsilon$}
              \State $\widetilde{w}=\frac{1}{2\pi\epsilon}\rm Arg\Big( \frac{\mathcal{F}^{sort}({\mathbf f}_{p, \epsilon})[h]}{\mathcal{F}^{sort}({\mathbf f}_{p, 0})[h]} \Big)$
              \State $a=\frac{1}{p}\mathcal{F}^{sort}({\mathbf f}_{p, 0})[h]$
              \State $R \gets R\cup ({\widetilde{w}}, a)$
              \EndIf
        \EndFor
        \State{prune small coefficients from $R$}
\State $i \gets i+1$
\EndWhile
\EndProcedure
\Statex
\end{algorithmic}
  \vspace{-0.4cm}
\end{algorithm}

\subsection{Multidimensioanl Problem Setting and Worst Case Scenario} \label{2.2}
 In this section, the multidimensional problem is introduced. Let us consider a function $f: {\mathbb R^d} \rightarrow {\mathbb C}$ such that
\begin{equation}\label{(9)}
	f({\mathbf t})~=~\sum_{j=1}^{k}a_j e^{2\pi i{{\mathbf w}_j}\cdot{{\mathbf t}}},
\end{equation}
where  ${\mathbf w}_j \in [-N/2, N/2)^d\cap{\mathbb Z}^d$ and $a_j\in{\mathbb C}$. That is, from (\ref{(1)}), $t$ is replaced by the $d$-dimensional phase or time vector ${\mathbf t}$, frequency $w_j$ is replaced by the frequency vector ${{\mathbf w}_j}$ and thus the operator between ${{\mathbf w}_j}$ and ${\mathbf t}$ is a dot product instead of simple scalar multiplication. We can see that this is a natural extension of the one-dimensional sparse problem. As in the 1D setting, if we find $a_j$ and ${\mathbf w}_j$, we recover the function $f$.

 However, since our time and frequency domain have changed, we cannot apply the previous algorithm directly. If we project the frequencies onto a line, then we can apply the former algorithm so that we can retain sublinear time complexity. Since the operator between frequency and time vectors is a dot product, we can convert the projection of frequencies to that of time. For example, we consider the projection onto the first axis, that is, we put the last $d-1$ elements of time vectors as $0$. If the projection is one-to-one, i.e., there is no collision, then we can apply the algorithm in Section \ref{2.1} to this projected function to recover the first element of frequency vectors. If there is a collision on the first axis, then we can try another projection onto $i$-th axis, $i=2, 3, \cdots, d$, until there are no collisions. We introduce in latter sections how to recover the corresponding remaining $d-1$ elements by extending the test to determine the occurrence of a collision in Section \ref{2.1}. Furthermore, to reduce the chance  of a collision through projections, we use an \lq\lq unwrapping method'' which unwraps frequencies onto a lower dimension guaranteeing a one-to-one projection.  There are both a \lq\lq full unwrapping'' and a \lq\lq partial unwrapping'' methods, which are explained in later sections.

 We shall call projections onto any one of the coordinate axes a {\em parallel projection}. The {\em worst case} is where there is a collision for every parallel projection. This obviously happens when a subset of frequency vectors form the vertices of a $d$-dimensional hypercube, but it can happen also with other configurations that require fewer vertices. Then our method cannot recover any of these frequency vectors via parallel projections. To resolve this problem, we introduce {\em tilted projections}: instead of simple projection onto axes we project frequency vectors onto tilted lines or planes so that there is no collision after the projection. We shall call this the {\em tilting method} and provide the details in the next section.  After introducing these projection methods, we explore which combination of these methods is likely to be optimal.

\section{TWO DIMENSIONAL SUBLINEAR SPARSE FOURIER ALGORITHM}\label{3}
\setcounter{equation}{0}

 As means of explanation, we introduce the two-dimensional case in this section and extend this to higher dimensions in Section \ref{4}.  The basic two-dimensional algorithm using a parallel projection is introduced in Section \ref{3.1}, the full unwrapping method is introduced in Section \ref{3.2} and the tilting method for the worst case scenarios is discussed in Section \ref{3.3}.

\begin{figure}[ht]
	\centering
  \includegraphics[width=0.9\textwidth, angle=0]{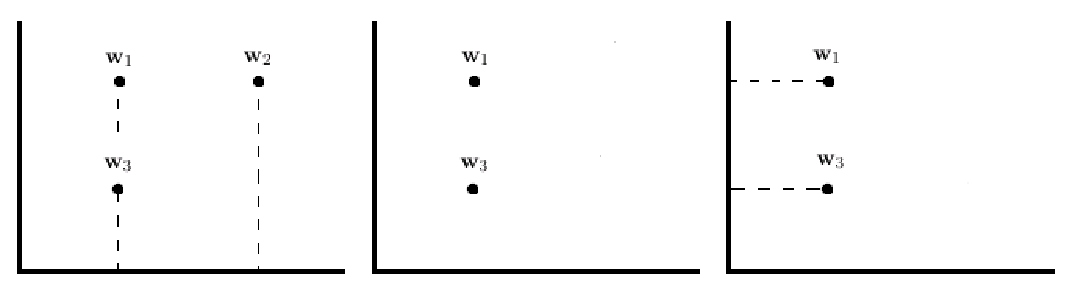}
	\caption{Process of the basic algorithm in 2D}
	\label{2Dbasic}
\end{figure}

\subsection{Basic Algorithm Using Parallel Projection} \label{3.1}
 Our basic two-dimensional sublinear algorithm excludes certain worst case scenarios. In most cases, we can recover frequencies in the 2-D plane by projecting them onto each horizontal axis or vertical axis.  Figure \ref{2Dbasic} is a simple illustration. Here we have three frequency vectors where $\mathbf w_1$ and $\mathbf w_3$ are colliding with each other when they are projected onto the horizontal axis, and $\mathbf w_1$ and $\mathbf w_2$ are when they are projected onto the vertical axis. The first step is to project the frequency vectors onto the horizontal axis and recover $\mathbf w_2$ and its corresponding coefficient $a_2$ only, since it is not colliding. By using two additional samples sets at the locations shifted along both axes of physical domain, each component of $\mathbf w_2$ can be recovered. After subtracting $\mathbf w_2$ from the data, we project the remaining frequency vectors onto the vertical axis and then find both $\mathbf w_1$ and $\mathbf w_3$.

 Now let us consider the generalized two-dimensional basic algorithm. Assume that we have a two-dimensional function $f$ with sparsity $k$ :
\begin{equation}\label{(10)}
	f({{\mathbf t}})~=~\sum^{k}_{j=1}a_j e^{2\pi i {{\mathbf w}_j}\cdot{{\mathbf t}}},\quad a_j\in {\mathbb C},\quad {{\mathbf w}_j} \in  
	\Big[ -\frac{N}{2} , \frac{N}{2} \Big)^2\cap \mathbb{Z}^2.
\end{equation}
For now, let us focus on one frequency vector with index $j'$ which is not collided with any other pairs when they are projected onto the horizontal axis. To clarify put ${{\mathbf t}}=(t_1, 0)$ with ${{\mathbf w}_j}=(w_{j1}, w_{j2})$ into (\ref{(10)}),
\begin{equation}\label{(11)}
	f^1(t_1)~:=~f(t_1,0)~=~\sum^{k}_{j=1}a_j e^{2\pi i {w_{j1}} {t_1}},
\end{equation}
which gives the same effect of parallel projection of frequency vectors.
Now, we can consider this function as a one-dimensional function $f^1$ so that we can use the original one dimensional sparse Fourier algorithm to find the first component of ${{\mathbf w}_{j'}}$. We get the samples ${\mathbf f}^1_{p,0}$ and ${\mathbf f}^1_{p,\epsilon}$ with and without shift by $\epsilon$. We can find these in the form of sequences in (\ref{(2)}), apply the DFT to them, and then recover the first component of the frequency pair and its coefficient as follows,
\begin{align}
	w_{j'1} &~=~ \frac{1}{2\pi \epsilon}\rm Arg\Big( \frac{{\mathcal F}({\mathbf f}^1_{p,\epsilon})[h]} {{\mathcal F}({\mathbf f}^1_{p,0})[h]} \Big), \nonumber\\
	a_{j'} &~=~ \frac{1}{p}{\mathcal F}({\mathbf f}^1_{p,0})[h].\label{(12)}
\end{align}
	By a projection onto the horizontal axis one can compute only the first component of $\mathbf{w}_{j'}$ and therefore one cannot simply subtract before projection onto the vertical axis. This is because the function $f^1$ is the linear combination of Fourier basis consisting only of the first components of each frequency pair.
	%Please give a more detailed explanation here.
	At the same time, we need to find the second component. In (\ref{(11)}), we replace $0$ by $\epsilon$. Then 
\begin{align}
	f^2(t_1)~ &:= ~ f(t_1,\epsilon) ~ = ~ \sum^{k}_{j=1}a_j e^{2\pi i ({w_{j1}}t_1+{w_{j2}}\epsilon)}\nonumber,\\
	{\mathcal F}({\mathbf f}^2_{p,\epsilon})[h]~ &= ~ pa_{j'}e^{2\pi i w_{j'2}\epsilon}\nonumber,\\
	w_{j'2}~ &= ~ \frac{1}{2\pi\epsilon} {\rm Arg}\Big( \frac{{\mathcal F}({\mathbf f}^2_{p,\epsilon})[h]}{ {\mathcal F}({\mathbf f}^1_{p,0})[h]} \Big),\label{(16)}
\end{align}
where ${\mathbf f}^2_{p,\epsilon}$ are samples shifted by $\epsilon$ in the vertical sense with rate $1/p$ from the function $f^2$. (\ref{(12)}) holds only when $w_{j'1}$ is the only one congruent to $h$ modulo $p$ among every first component of $k$ frequency pairs  and (\ref{(16)}) holds only when the previous condition is satisfied and ${\mathbf w}_{j'}=(w_{j'1},w_{j'2})$ does not collide with other frequency pairs from the parallel projection.

 Now we have two kinds of possible collisions. The first one is from taking modulo $p$ after the parallel  projection and the second one is from the projection. Thus we need two tests. To determine whether there are both kinds of collisions, we use similar tests as (\ref{(8)}). If there are at least two different $w_{j1}$ congruent to $h$ modulo $p$, then the second equality in the following is not satisfied for almost all $\epsilon$, just as (\ref{(8)}),
\begin{equation}\label{(17)}
	\frac{\vert{\mathcal F}({\mathbf f}^1_{p,\epsilon})[h]\vert}{\vert{\mathcal F}({\mathbf f}^1_{p,0})[h]\vert}~=~ \frac{\vert p \sum_{w_{j1} = h(\bmod p)} a_j e^{2 \pi i \epsilon w_{j1}} \vert}{\vert p \sum_{w_{j1} = h(\bmod p)} a_j  \vert} ~=~ 1.
\end{equation}
Second, if there is a collision from the projection, i.e., the first components $w_{j1}$'s of at least two frequency vectors are identical and the corresponding $w_{j2}$'s are different, the following second equality  does not hold for almost all $\epsilon$,
\begin{equation}\label{(18)}
	\frac{\vert{\mathcal F}({\mathbf f}^2_{p,\epsilon})[h]\vert}{\vert{\mathcal F}({\mathbf f}^1_{p,0})[h]\vert}~=~ \frac{\vert p \sum_{w_{j1} = h (\bmod p)} a_j e^{2 \pi i \epsilon w_{j2}} \vert}{\vert p \sum_{w_{j1} = h (\bmod p)} a_j  \vert} ~=~ 1.
\end{equation}
The two tests above are both satisfied only when there is no collision both from taking modulo $p$ and the projection. We use these for the complete recovery of the objective frequencies.

 So far we project the frequencies onto the horizontal axis. After we find the non-collided frequencies from the first projection, we subtract a function consisting of found frequencies and their coefficients from the original function $f$ to get a new function. Next we project this new function onto the vertical axis and do a similar process. The difference is to exchange $1$ and $2$ in the super-indices and sub-indices respectively in (\ref{(11)}) through (\ref{(18)}). Again, find the remaining non-collided frequencies, change the axis again and keep doing this until we recover all of the frequencies.

\subsection{Full Unwrapping Method} \label{3.2}
 We introduce another kind of projection which is one-to-one. The full unwrapping method uses one-to-one projections onto one-dimensional lines instead of the parallel projection onto axes from the previous method. We consider the $k$ pairs of frequencies $(w_{j1}, w_{j2})$, $j=1,2,\cdots,k$ and transform them as follows
\begin{equation}\label{(19)}
	(w_{j1}, w_{j2}) ~\rightarrow~ w_{j1}+Nw_{j2}.
\end{equation}
This transformation in frequency space can be considered as the transformation in phase or time space. That is, from the function in (\ref{(10)})
\begin{equation}\label{(20)}
	g(t)~:=~f(t, Nt)~=~\sum^{k}_{j=1}a_j e^{2\pi i (w_{j1}+Nw_{j2})t}.
\end{equation}
The function $g(t)$ is a one-dimensional function with sparsity $k$ and bandwidth bounded by $N^2$. We can apply the algorithm in Section \ref{2.1} on $g$ so that we recover $k$ frequencies of the form on the right side of the arrow in (\ref{(19)}). Whether unwrapped or not, the coefficients are the same, so we can find them easily. In the end we need to wrap the unwrapped frequencies to get the original pairs. Remember that unwrapping transformation is one-to-one. Thus we can wrap them without any collisions.

 Since the pairs of the frequencies are projected onto the one-dimensional line directly, we call this method the ``full unwrapping method". Problem with this method occurs when the  dimension $d$ gets large. From the above description, we see that after the one-to-one unwrapping the total bandwidth of the two dimensional signal increases from $N$ in each dimension to $N^2$. If the full unwrapping method is applied to a function in $d$-dimensions, then to guarantee the one-to-one transformation the bandwidth will be $N^{d}$.  Theoretically this does not matter. However, since $\epsilon$ is dependent on the bandwidth, in the case where $d$ is large, we need to consider the limit of machine precision for practical implementations.  As a result, we need to introduce the partial unwrapping method to prevent the bandwidth from becoming too large. The partial unwrapping method is discussed in Section \ref{4}.
\begin{figure}[ht]
	\centering
  \includegraphics[width=0.9\textwidth, angle=0]{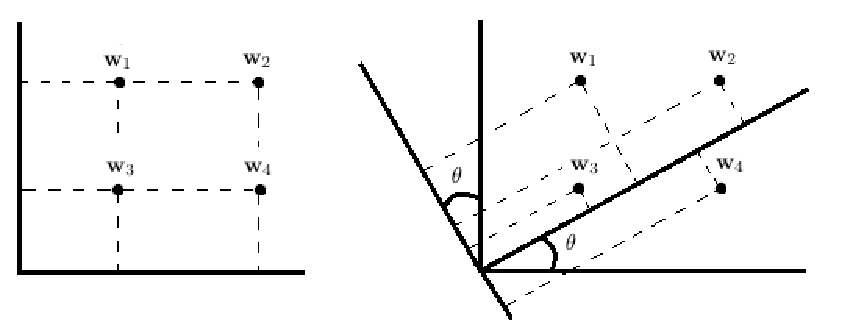}
	\caption{Worst case scenario in $2D$ and solving it through tilting}
	\label{2Dtilting}
\end{figure}

\subsection{Tilting Method for the Worst Case} \label{3.3}

\vspace{-0.2cm}
\alglanguage{pseudocode}
\begin{algorithm}[H]
	\small
	\caption{2D Sparse Fourier Algorithm with Tilting Method Pseudo Code}
	\label{Algorithm 3}
	\begin{algorithmic}[1]
		\Procedure{$\mathbf{2DTiltedPhaseshift}$}{}\\
		{\textbf{Input:}}{$f, c, k, N, d, \epsilon$, integers $base$, $height$, $hypo$}\\
		{\textbf{Output:}}{$R$}
		\State $R \gets \emptyset$
		\State $i \gets 1$
		\State $\cos \gets {base}$, $\sin \gets {height} $
		\While { $|R|<k$}
		\State $k^{\ast} \gets k - |R|$
		\State $p \gets {\it{i}\text{-th}}$ prime number $\geq c k^{\ast}$
		\State $ m \gets (i$ mod 2)+1
		\State $g({\mathbf{t}})=\sum_{({\mathbf w},a_{\mathbf w})\in R} a_{\mathbf w} e^{2 \pi i {\mathbf w}\cdot {\mathbf t}}$
		\For {$n = 1 \to 2$}
		\For {$h = 1 \to p$}
		\State $m' \gets m$ mod $2$, $m'' \gets m+1$ mod $2$, $n' \gets n$ mod $2$, $n'' \gets n+1$ mod $2$
		\State ${f}^{m,n}_{p,\epsilon}[h] = $
		\algrenewcommand\algorithmicindent{4.5em}%
		\Statex \hskip \algorithmicindent $f((\frac{h-1}{p}m'+\epsilon n') \cos+ (\frac{h-1}{p}m''+\epsilon n'') \sin, -(\frac{h-1}{p}m'+\epsilon n') \sin+ (\frac{h-1}{p}m''+\epsilon n'') \cos)$
		\Statex \hskip \algorithmicindent $- g(\frac{h-1}{p}{\mathbf e}_{m}+\epsilon{\mathbf e}_{n} )$
		\algrenewcommand\algorithmicindent{1.0em}%
		\State ${f}^{m,n}_{p,0}[h] = f(\frac{h-1}{p}m' \cos +\frac{h-1}{p}m'' \sin , -\frac{h-1}{p}m' \sin +\frac{h-1}{p}m'' \cos ) - g(\frac{h-1}{p}{\mathbf e}_m)$
		\EndFor
		\State $\mathcal{F}({f}^{m,n}_{p, \epsilon})=FFT({f}^{m,n}_{p, \epsilon})$
		\State $\mathcal{F}({f}^{m,n}_{p, 0})=FFT({f}^{m,n}_{p, 0})$
		\State $\mathcal{F}^{sort}({f}^{m,n}_{p, 0})=SORT(\mathcal{F}({f}^{m,n}_{p, 0}))$
		\EndFor
		\For {$h = 1 \to k^*$}
		\State $\ell \gets 0$
		\For {$n = 1 \to 2$}
		\If {$\Big|\frac{|\mathcal{F}^{sort}({f}^{m,n}_{p, 0})[h]|}{|\mathcal{F}^{sort}({f}^{m,n}_{p, \epsilon})[h]|}-1\Big|<\epsilon$}
		\State $\ell \gets \ell +1$
		\EndIf
		\State $\tilde{w}_{n}=\frac{1}{2\pi\epsilon}{\rm Arg}\Big( \frac{\mathcal{F}^{sort}({f}^{m,n}_{p, \epsilon})[h]}{\mathcal{F}^{sort}({f}^{m,n}_{p, 0})[h]} \Big)$
		\State $a=\frac{1}{p}\mathcal{F}^{sort}({f}^{m,n}_{p, 0})[h]$
		\EndFor
		\If {$\ell==2$}
		\State $R \gets R\cup ({\mathbf{ \tilde{w}}}, a)$
		\EndIf
		\EndFor
		\State{prune small coefficients from $R$}
		\State $i \gets i+1$
		\EndWhile
		\State $\cos \gets \frac{base}{hypo}$, $\sin \gets \frac{height}{hypo} $
		\State {rotate each ${\mathbf{ \tilde{w}}}$ back to ${\mathbf {w}}$ using a matrix [$\cos$ $\sin$;${-\sin}$ $\cos$] and restore it in $R$}
		\EndProcedure
		\Statex
	\end{algorithmic}
	\vspace{-0.4cm}
\end{algorithm}
 Up till now, we have assumed that we do not encounter the worst case, i.e., that we do not encounter the case where any frequency pair has collisions from the parallel projection for all coordinate axes. This makes the algorithm break down. The following method is for finding those frequency pairs. Basically, we rotate axes of the frequency plane and thus use a projection onto a one-dimension system which is a tilted line with the tilt chosen so that there are no collisions. If the horizontal and vertical axes are rotated with angle $\theta$ then the frequency pair ${\mathbf w}_j=(w_{j1}, w_{j2})$ can be relabeled with new coordinates as the right side of the following

\begin{equation}\label{(21)}
	(w_{j1}, w_{j2}) ~\rightarrow~ (w_{j1}\cos\theta - w_{j2} \sin\theta , w_{j1} \sin\theta + w_{j2} \cos\theta ).
\end{equation}
In phase-sense, this rotation can be written as
\begin{align}
	g({\tilde t}_1, {\tilde t}_2)& ~:=~
f({\tilde t}_1 \cos\theta + {\tilde t}_2 \sin\theta, - {\tilde t}_1 \sin\theta + {\tilde t}_2 \cos \theta)\nonumber \\
	&~=~\sum^{k}_{j=1}a_j e^{2\pi i\{ w_{j1}( {\tilde t}_1 \cos\theta + {\tilde t}_2 \sin\theta)+w_{j2}(- {\tilde t}_1 \sin\theta + {\tilde t}_2 \cos \theta)\}},\nonumber\\
	&~=~\sum^{k}_{j=1}a_j e^{2\pi i \{(w_{j1} \cos\theta - w_{j2}\sin\theta ){\tilde t}_1+(w_{j1} \sin\theta + w_{j2} \cos\theta ){\tilde t}_2\}}.\label{(22)}
\end{align}
We can apply the basic algorithm in Section \ref{3.1} to the function $g$ to get the frequency pairs in the form of the right side of the arrow in (\ref{(21)}).

 One problem we face is that the components of the projected frequency pairs should be integers to apply the method, since we assume the integer frequencies in the first place. Any irrational $\tan\theta$ makes the projected frequencies become irrational, which is undesirable although it guarantees the injectivity. Thus, we should try a rational $\tan\theta$, and to make them integers. It is unavoidable that this process of rotation increases the bandwidth. To minimize this effect we choose to multiply by the least common multiple of the denominators of $\sin\theta$ and $\cos\theta$. We assume the following,  
\begin{equation}\label{(25)}
	\sin\theta ~=~ \frac{a}{c},\quad \cos\theta ~=~ \frac{b}{c},\quad \gcd(a,c)~=~\gcd(b,c)~=~1,
\end{equation}
where $a$, $b$ and $c$ are integers.
Multiplying $c$ to both inputs in the right-hand side of (\ref{(22)}) we obtain
\begin{align}
	{\hat g}({\tilde t}_1, {\tilde t}_2)&~:=~f(c({\tilde t}_1 \cos\theta + {\tilde t}_2 \sin\theta), c(- {\tilde t}_1  \sin\theta +{\tilde t}_2\cos \theta))\nonumber\\
	&~=~\sum^{k}_{j=1}a_j e^{2\pi i \{(cw_{j1} \cos\theta -c w_{j2}\sin\theta ){\tilde t}_1+(c w_{j1} \sin\theta +cw_{j2} \cos\theta ){\tilde t}_2\}}.\label{(27)}
\end{align}

As long as there is no collision for at least one projection, the frequency pairs, $(c w_{j1} \cos\theta - c w_{j2} \sin \theta , cw_{j1} \sin \theta + c w_{j2} \cos \theta )$, can be found by applying the basic algorithm in Section \ref{3.2} on ${\hat g}$. Due to the machine precision the integer $c$ should not be too large, or the bandwidth gets too large resulting in the failure of the algorithm. If four pairs of frequencies are at the vertices of a rectangle aligned with coordinate axes before the rotation, then they are not aligned after the rotation with $0<\theta<\pi/2$. Thus we can assure finding the frequencies whether they are in the worst case or not.

The pseudo code of the 2D tilting method is shown in Algorithm \ref{Algorithm 3}. The lines 14 and 15 mean that each frequency pair $(w_{j1},w_{j2})$ is rotated by a matrix [$\cos$ ${-\sin}$; $\sin$ $\cos$] and scaled to make the rotated components integers. Thus we first find the frequency pairs in the form of ${\mathbf {\tilde{w}}}=( w_{j1}\cos-w_{j2}\sin , w_{j1}\sin + w_{j2}\cos )$ and after finding all of them, we rotate them back into the original pairs with the matrix [$\cos$ $\sin$; ${-\sin}$ $\cos$] in line 39.

This tilting method is a straight forward way to resolve the worst case problem. First, we recover the frequencies as much as possible from the basic parallel projection method. If we cannot get any frequency pairs for several projections switching among each axis then, assuming that the worst case happens, we apply the tilting method with an angle so that all remaining frequency pairs are found. We only introduced the tilting method in the two-dimensional case, but the idea of rotating the axes can be extended to the general $d$-dimensional case with some effort. On the other hand, we may notice that the probability of this worst case is very low, especially when the number of dimensions $d$ is large. Thus, as we recover the frequencies as much as possible from the basic algorithm.  Its details are shown in Section \ref{4}.

\section{PARTIAL UNWRAPPING METHOD FOR HIGH DIMENSIONAL ALGORITHM}
\label{4}
\setcounter{equation}{0}

 In this section we present the {\em partial unwrapping} method for a sublinear sparse Fourier algorithm for very high dimensional data. As we have already mentioned, while full unwrapping converts a multi-dimensional problem into a single dimensional problem, it is severely limited in its viability when the dimension is large or when the bandwidth is already high because of the increased bandwidth. Partial unwrapping is introduced here to  overcome this problem and other problems.  In Section \ref{4.1} we give a four dimensional version of the algorithm using the partial unwrapping method as well as a generalize it to $d$ dimension.  In Section \ref{4.2}, the probability of the worst case in $d$ dimension is analyzed.

\subsection{Partial Unwrapping Method} \label{4.1}
 To see the benefit of partial unwrapping we need to examine the main difficulties we may encounter in developing sublinear sparse Fourier algorithms. For this let us consider a hypothetical case of sparse FFT where we have $k=100$ frequencies in a 20-dimensional Fourier series distributed in $[-10, 10)^{20}$. When we perform the parallel projection method, because the bandwidth is small, there will be a lot of collisions after the projections. It is often impossible to separate any frequency after each projection, and the task could thus not be completed. This, ironically, is a {\em curse of small bandwidth} for sparse Fourier algorithm. On the other hand, if we do the full unwrapping we would have increased the bandwidth to $N=20^{20}$, which is impossible to do within reasonable accuracy because $N$ is too large.

However, a partial unwrapping would reap the benefit of both worlds. Let us now break down the 20 dimensions into 5 lower 4-dimensional subspaces, namely we write
$$
    [-10,10)^{20} = \left( [-10,10)^4\right)^{5}.
$$
In each subspace we perform the full unwrapping, which yields bandwidth $N=20^4=160,000$ in the subspace. This bandwidth $N$ is large enough compared with $k$, so when projection method is used there is a very good probability that collision will occur only for a small percentage of the frequencies, allowing them to be reconstructed. On the other hand, $N$ is not so large that the phase-shift method will incur significant error.

One of the greatest advantage of partial unwrapping is to turn the curse of dimensionality into the {\em blessing} of dimensionality.

Note that in the above example, the 4 dimensions that for any of the subspaces do not have to follow the natural order. By randomizing (if necessary) the order of the dimensions it may achieve the same goal as the tilting method would.
Also note that the dimension for each subspace needs not be uniform. For example, we can break down the above 20-dimensional example into  four $3$-dimensional subspaces and two $4$-dimensional subspaces, i.e.
$$
     [-10,10)^{20} = \left( [-10,10)^3\right)^{4} \times \left( [-10,10)^4\right)^{2}.
$$
This will lead to further flexibility.

\subsubsection{Example of 4-D Case}\label{4.1.1}
 Before introducing the generalized partial unwrapping algorithm for dimension $d$, let us think about the simple case of $4$ dimensions. We assume that $k$ frequency vectors are in $4$-dimensional space ($d=4$). Then, a function $f$ constructed from these frequency vectors is as follows,
\begin{equation}\label{(28)}
	f({\mathbf t})=\sum^{k}_{j=1}a_j e^{2\pi i {\mathbf w}_j\cdot{\mathbf t}},\ a_j\in {\mathbb C},\ {\mathbf w}_j \in  {\Big{(}\Big[-\frac{N}{2},\frac{N}{2}{\Big )}\cap{\mathbb Z}\Big{)}}^4.
\end{equation}
Since $4=2\times 2$, the frequency pairs of the two-two dimensional spaces are both unwrapped onto one-dimensional spaces. Here, $4$ dimensions is projected onto $2$ dimensions as follows
\begin{align}
	g(t_1,t_2)&~:=~f(t_1, Nt_1, t_2, Nt_2)\nonumber\\
	&~=~\sum^{k}_{j=1}a_j e^{2\pi i\{ (w_{j1}+Nw_{j2})t_1+ (w_{j3}+Nw_{j4})t_2\}}\nonumber\\
	&~=~ \sum^{k}_{j=1}a_j e^{2\pi i(\tilde{w}_{j1}t_1+ \tilde{w}_{j2}t_2)},\label{(31)}
\end{align}
where $\tilde{w}_{j1}=w_{j1}+Nw_{j2}$ and $\tilde{w}_{j2}=w_{j3}+Nw_{j4}$. Note that this projection is one-to-one so as to guarantee the inverse transformation.

 Now we can apply the basic projection method in Section \ref{3.1} to this function $g$ re-defined as the $2$-dimensional one. To make this algorithm work, ${\mathbf {\tilde{w}}}_j=(\tilde{w}_{j1},\tilde{w}_{j2})$ should not collide with any other frequency pair after the projection onto either the horizontal or vertical axes. If not, we can consider using the tilting method. After finding all the frequencies in the form of $(\tilde{w}_{j1},\tilde{w}_{j2})$, it can be transformed to $(w_{j1},w_{j2},w_{j3},w_{j4})$.

\subsubsection{Generalization}\label{4.1.2}
 We introduce the final version of the multidimensional algorithm in this section. Its pseudo code and detailed explanation are given in Algorithm \ref{Algorithm2} and Section \ref{5.0}, respectively. We start with a $d$-dimensional function $f$,
\begin{equation}\label{(32)}
	f({\mathbf t})~=~\sum^{k}_{j=1}a_j e^{2\pi i {\mathbf w}_j\cdot{\mathbf t}},\quad  a_j\in {\mathbb C},\quad  {\mathbf w}_j \in  {\Big{(}\Big[-\frac{N}{2},\frac{N}{2}{\Big )}\cap{\mathbb Z}\Big{)}}^d.
\end{equation}
Let us assume that $d$ can be divided into $d_1$ and $d_2$ - the case of $d$ being a prime number will be mentioned at the end of this section. The domain of frequencies can be considered as  ${([-N/2,N/2)\cap{\mathbb Z})}^d=({([-N/2,N/2)\cap{\mathbb Z})}^{d_1})^{d_2}$ and ${([-N/2,N/2)\cap{\mathbb Z})}^{d_1}$ will be reduced to one dimension, as $d_1$ is in the 4 dimensional case. Each of the $d_1$ elements of a frequency vector, ${\mathbf w}_j =  (w_{j1}, w_{j2}, \cdots, w_{jd} )$, is unwrapped as
\begin{align}
	& (w_{j(d_1q+1)}, w_{j(d_1q+2)}, w_{j(d_1q+3)},  \cdots, w_{j(d_1q+d_1)})\nonumber\\
	&~\rightarrow~ w_{j(d_1q+1)}+Nw_{j(d_1q+2)}+N^2w_{j(d_1q+3)}+\cdots+N^{d_1-1}w_{j(d_1q+d_1)}\nonumber\\
	&~=:~\tilde{w}_{j(q+1)}\label{(35)}
\end{align}
with $q=0,1,2,\cdots,d_2-1$, increasing the respective bandwidth from $N$ to $N^{d_1}$ and having injectivity. We rewrite this transformation in terms of the phase. With ${\mathbf t}=(t_1,t_2, \cdots, t_d)$ and put the following into $t_\ell$
\begin{equation}\label{(36)}
	N^{R(\ell-1,d_1)}\tilde{t}_{Q(\ell,d_1)}
\end{equation}
for all $\ell = 1, 2, \cdots, d$, where $R(\ell-1,d_1)$ and $Q(\ell,d_1)$ are the remainder  from dividing $\ell-1$ by $d_1$ and quotient from dividing $\ell$ by $d_1$  respectively, and ${\mathbf {\tilde{t}}}=(\tilde{t}_1, \tilde{t}_2,\cdots, \tilde{t}_{d_2})$ is a phase vector in $d_2$ dimensions after projection. Define a function $g$ on $d_2$ dimension as
\begin{align}
	g({\mathbf {\tilde{t}}})&~:=~ f(\cdots, N^{R(\ell-1,d_1)}\tilde{t}_{Q(\ell,d_1)}, \cdots)\nonumber\\
	&~=~\sum^{k}_{j=1}a_j e^{2\pi i \sum^{d_2-1}_{q=0} \big( \sum^{d_1-1}_{r =0} w_{j (d_1q+r+1)}N^{r}\big)\tilde{t}_{q+1}},\label{(38)}
\end{align}
where $N^{R(\ell-1,d_1)}\tilde{t}_{Q(\ell,d_1)}$ is the $\ell$th element of the input of $f$. If we project frequency vectors of $g$ onto the $m$th axis, then the $n$th element of a frequency vector ${\mathbf {\tilde{w}}}_j$ can be found in the following computation,
\begin{align}
	{\mathbf g}^{m,n}_{p, \epsilon} &~=~ \left( g(0{\mathbf e}_m + \epsilon{\mathbf e}_n), g\left(\frac{1}{p}{\mathbf e}_m + \epsilon{\mathbf e}_n\right), \cdots, g\left(\frac{p-1}{p}{\mathbf e}_m + \epsilon{\mathbf e}_n\right)\right) \nonumber\\
	\tilde{w}_{jn} &~=~ \frac{1}{2\pi \epsilon}{\rm Arg}\Big( \frac{{\mathcal F}({\mathbf g}^{m,n}_{p,\epsilon})[h]} {{\mathcal F}({\mathbf g}^{m,n}_{p,0})[h]} \Big) \nonumber\\
	a_j &~=~ \frac{1}{p}{\mathcal F}({\mathbf g}^{m,n}_{p,0})[h],\label{(41)}
\end{align}
where ${\mathbf e}_m$ is the $m$-th unit vector with length $d_2$, i.e., all elements are zero except the $m$-th one with entry $1$.  (\ref{(41)}) holds as long as $\tilde{w}_{jn}$ is the only one congruent to $h$ modulo $p$ among all $n$-th elements of the frequency vectors and ${\mathbf {\tilde{w}}}_j$ does not collide with any other frequency vector due to the projection onto the $m$-th axis. The test for checking whether these conditions are satisfied is
\begin{equation}\label{(42)}
	\frac{\vert{\mathcal F}({\mathbf g}^{m,n}_{p,\epsilon})[h]\vert}{\vert{\mathcal F}({\mathbf g}^{m,n}_{p,0})[h]\vert} ~=~1
\end{equation}
for all $1\leq n \leq d_2$. The projections onto the $m$-th axis, where $m=1, \cdots,d_2$, take turns until we recover all frequency vectors and their coefficients. After that we wrap the unwrapped frequency vectors up from $d_2$ to $d$ dimension. Since the unwrapping transformation is one-to-one, this inverse transformation is well-defined.

  So far, we assumed that dimension $d$ can be divided into two integers, $d_1$ and $d_2$. For the case that $d$ is a  prime number or both $d_1$ and $d_2$ are so large that the unwrapped data has a  bandwidth such that $\epsilon$ is below the machine precision, a strategy of divide and conquer can be applied. In that case we can think about applying partial unwrapping method in a way that each unwrapped component has a different size of bandwidth. If $d$ is $3$, for example, then we can unwrap the first two components of the frequency vector onto one dimension and the last one lies in the same dimension. In that case, the unwrapped data is in two dimensions, and the bandwidth of the first component is bounded by $N^2$ and that of second component is bounded by $N$. In this case we can choose a shift $\epsilon< 1/N^2$ where $N^2$ is the largest bandwidth. We can extend this to the general case, so the partial unwrapping method has a variety of choices balancing the bandwidth and machine precision.

\subsection{Probability of Worst Case Scenario} \label{4.2}
In this section, we give an upper bound of the probability of the worst case assuming that we randomly choose a partial unwrapping method. As addressed in the Section \ref{4.1}, there is flexibility in choosing certain partial unwrapping method. Assuming a certain partial unwrapping method and considering a stronger condition to avoid its failure, we can find the upper bound of the probability of the worst case where there is a collision for each parallel projection.

For simple explanation, consider a two dimensional problem. Choosing the first frequency vector $(w_{11},w_{12})$ on a two dimensional plane, if the second frequency vector, $(w_{21},w_{22})$, is not on the vertical line crossing $(w_{11},0)$ and the horizontal line crossing $(0,w_{12})$, then the projection method works. Then if the third frequency vector is not on four lines, those two lines mentioned before, the vertical line crossing $(w_{21},0)$ and the horizontal line crossing $(0,w_{22})$, then again the projection method works. We keep choosing next frequency vector in this way, excluding the lines containing previous frequencies. Thus, letting such event $A$, the probability that the projection method fails is bounded above by $1-\mathbb{P}(A)$.

Generally, let us assume that we randomly choose a partial unwrapping, without loss of generality, the total dimension is $d=d_1+d_2+\cdots +d_r$ where $r$ is the number of subspaces and $d_1, d_2,\cdots, d_r$ are the dimensions of each subspace.
That is, partially unwrapped frequency vectors are in $r<d$ dimension and each bandwidth is $N^{d_1}, N^{d_2},\cdots, N^{d_r}$, respectively, which is integer strictly larger than $1$. Then, the failure probability of projection method is bounded above by 
\begin{align}
	1-\prod_{j=1}^{k} \mathbb{P}(A_j)& ~\leq~ 1-\prod_{j=1}^{k} \frac{{N^d-(j-1)(N^{d_1}+N^{d_2}+\cdots+N^{d_r})}}{N^d}\nonumber\\
	&~=~1- \frac{1}{\tau^{k}}\frac{\tau !}{(\tau-k)!} \quad \left(\tau:=\frac{N^d}{N^{d_1}+N^{d_2}+\cdots+N^{d_r}}\right)\nonumber\\
	&~\sim~ 1- \frac{1}{\tau^k} \sqrt{\frac{\tau}{\tau-k}}\frac{\left( \frac{\tau}{e} \right)^{\tau}}{\left( \frac{\tau-k}{e}\right)^{\tau-k}} \quad(\text{by Sterling's formula})\nonumber\\
	&~=~1-\frac{1}{e^k}\left( 1-\frac{k}{\tau}\right)^{-\frac{\tau}{k} \cdot k} \left( 1-\frac{k}{\tau} \right)^{k-\frac{1}{2}}
\end{align}
where $A_j$ is the event that we choose $j$th frequency not on the lines, crossing formerly chosen frequency vectors and parallel to each coordinate axis. Noting $N^d=N^{d_1}\times{N^{d_2}}\times \cdots \times {N^{d_r}}$, sparsity $k$ is relatively small compared to $N^d$, and $\tau$ is large, we can see that the upper bound above gets closer to $0$ as $d$ or $N$ grows to infinity.

\section{Analysis}\label{sec:analysis}
\setcounter{equation}{0}
In this section, we analyze the performance of our algorithms suggested. We will prove that the tilting method works well in two dimensions but explain that it is not simple to extend the idea to the general high dimensional setting. However, it was shown in Section \ref{4.2} that the probability of the worst-case scenario is extremely small. Furthermore, the average-case runtime and sampling complexity is shown in this section under the assumption that there are no frequency vectors forming a worst-case scenario.
In \cite{Plonka2013}, it is conjectured that there exist some angles $\theta$ which suffice to recover all two-dimensional frequency vectors. 
In the following theorem, we identify a class of angles with which the tilting method in 2D recovers all frequency pairs even though they form a worst-case scenario. 
%
%\todo[inline]{Review 7. Section 4.1: Is it possible to compute an adapted tilted line in order to nd out which candidates of frequency pairs obtained from horizontal projections are indeed occurring? A similar idea has been used in the bivariate case in G.Plonka, M. Wischerho.How many Fourier samples are needed for real function reconstruction? Journal of Applied Mathematics and Computing 42 (2013), 117-137, see Section 4.2. However, also here the problem of neglecting sums can occur.}

\begin{theorem}
	Let ${\mathbf w}_j=(w_{j1},w_{j2}) \subset \left[-\frac{N}{2},\frac{N}{2} \right)^2 \cap \mathbb{Z}^2 $ for $j\in \{1,2,\cdots,k \}$. If $\tan \theta=\frac{a}{b}$ such that $c>b>a$ are Pythagorean triples where $b>2N$ and $a$ are relative primes, then all $(cw_{j1}\cos\theta -cw_{j2}\sin\theta,cw_{j1}\sin\theta +cw_{j2}\cos \theta)$ rotated by $\theta$ does not collide with any other pair through the parallel projection. Thus, all rotated pairs can be identified by the parallel projection method.
	\label{thm:tilting}
\end{theorem}

\vspace{-0.2cm}
\alglanguage{pseudocode}
\begin{algorithm}[H]
	\small
	\caption{Multidimensional Sparse Fourier Algorithm Pseudo Code}
	\label{Algorithm2}
	\begin{algorithmic}[1]
		\Procedure{$\mathbf{MultiPhaseshift}$}{}\\
		{\textbf{Input:}}{$f, c, k, N, d, d_1, d_2, \epsilon$}\\
		{\textbf{Output:}}{$R$}
		\State $R \gets \emptyset$
		\State $i \gets 1$
		\While { $|R|<k$}
		\State $k^{\ast} \gets k - |R|$
		\State $p \gets {\it{i}\text{-th}}$ prime number $\geq c k^{\ast}$
		\State $ m \gets (i$ mod $d_2)+1$
		\State $g({\mathbf{t}})=\sum_{({\mathbf w},a_{\mathbf w})\in R} a_{\mathbf w} e^{2 \pi i {\mathbf w}\cdot {\mathbf t}}$
		\For {$n = 1 \to d_2$}
		\For {$h = 1 \to p$}
		\State ${f}^{m,n}_{p,\epsilon}[h] = f(\sum^{d_1}_{\ell=1} N^\ell \frac{h-1}{p}{\mathbf e}_{d_1(m-1)+\ell}+ \epsilon \sum^{d_1}_{\ell=1} N^\ell {\mathbf e}_{d_1(n-1)+\ell} ) - g(\frac{h-1}{p}{\mathbf e}_m+\epsilon{\mathbf e}_n)  $
		\State ${f}^{m,n}_{p,0}[h] = f(\sum^{d_1}_{\ell=1} N^\ell \frac{h-1}{p}{\mathbf e}_{d_1(m-1)+\ell} ) - g(\frac{h-1}{p}{\mathbf e}_m)$
		\EndFor
		\State $\mathcal{F}({f}^{m,n}_{p, \epsilon})=FFT({f}^{m,n}_{p, \epsilon})$
		\State $\mathcal{F}({f}^{m,n}_{p, 0})=FFT({f}^{m,n}_{p, 0})$
		\State $\mathcal{F}^{sort}({f}^{m,n}_{p, 0})=SORT(\mathcal{F}({f}^{m,n}_{p, 0}))$
		\EndFor
		\For {$h = 1 \to k^*$}
		\State $\ell \gets 0$
		\For {$n = 1 \to d_2$}
		\If {$\Big|\frac{|\mathcal{F}^{sort}({f}^{m,n}_{p, 0})[h]|}{|\mathcal{F}^{sort}({f}^{m,n}_{p, \epsilon})[h]|}-1\Big|<\epsilon$}
		\State $\ell \gets \ell +1$
		\EndIf
		\State $\tilde{w}_{n}=\frac{1}{2\pi\epsilon}{\rm Arg}\Big( \frac{\mathcal{F}^{sort}({f}^{m,n}_{p, \epsilon})[h]}{\mathcal{F}^{sort}({f}^{m,n}_{p, 0})[h]} \Big)$
		\State $a=\frac{1}{p}\mathcal{F}^{sort}({f}^{m,n}_{p, 0})[h]$
		\EndFor
		\If {$\ell==d_2$}
		\State $R \gets R\cup ({\mathbf {\tilde{w}}}, a)$
		\EndIf
		\EndFor
		\State{prune small coefficients from $R$}
		\State $i \gets i+1$
		\EndWhile
		\State {inverse-transform each ${\mathbf {\tilde{w}}}$ in $d_2$-D to ${\mathbf {w}}$ $d$-D} and restore it in $R$
		\EndProcedure
		\Statex
	\end{algorithmic}
	\vspace{-0.4cm}
\end{algorithm}

\begin{proof}
	Suppose that any two frequency pairs  ${\mathbf w}_j=(w_{j1},w_{j2})$  and  ${\mathbf w}_{j'}=(w_{j'1},w_{j'2})$  cannot be recovered through the tilting method. This implies that the slope of the line crossing ${\mathbf w}_j$ and ${\mathbf w}_{j'}$ is perpendicular to $\tan \theta$, and thus those two pairs collide at least once if they are projected onto each principal axis rotated by the angle $\theta$. This results in the fact that $\frac{w_{j2}-w_{j'2}}{w_{j1}-w_{j'1}}$ is either $\frac{a}{b}$ or $-\frac{b}{a}$. However, this is a contradiction since $-N < w_{j1}-w_{j'1}, ~ w_{j2}-w_{j'2} < N$, and $a$ and $b>2N$ are relative primes.
\end{proof}
Theorem \ref{thm:tilting} implies that all frequency pairs can be distinguished by the tilting method with only one proper angle $\theta$. It is natural to think about extending this idea to the high-dimensional setting. However, it is not as easy to find a proper slope as in two dimensions using the Pythagorean triples. Moreover, even though we consider choosing a finite set of random angles guaranteeing that it includes the proper one, each line crossing between two arbitrary vectors in general $D$ dimensions has infinitely many lines perpendicular to it so that such a finite set is difficult to find again. 
It should be noted that there is a variety of the worst-case scenario where the relatively simple tilting method still works. For example, if the frequency vectors are located on the two dimensional subspace, then the tilting method applied on the subspace would recover all the frequencies. There still exist trickier worst-case scenarios, however, such that the frequencies form a $D$-dimensional cube or more complicated structures. It is still worth exploring further about the tilting method.
%-------if we get the proper worst case analysis
On the other hand, we also have shown in Section \ref{4.2} that the worst-case scenario rarely happens in extremely high dimensions. 

In order to obtain the average-case runtime and sampling complexity of the parallel projection method under the assumption that there is no worst-case scenario, we utilize the probability recurrence relation as in \cite{lawlor2013adaptive}. We remind readers that each entry of the frequency vectors should be distinguished modulo $N$ due to parallel projection as well as be distinguished modulo $p$. Since $c=5$ determining $p$ is shown in \cite{lawlor2013adaptive} to ensure that $90\%$ of frequencies are isolated modulo $p$ on average, at least $90\%$ of frequencies are isolated modulo $N$ on average if $N\geq p$. Taking the union bound of failure probabilities of isolation modulo $p$ and $N$ yields at least $80\%$ of frequencies are isolated both modulo $N$ and $p$ at the same time. Algorithm \ref{Algorithm2} with $d_1=1$ and $d_2=d$ implements the direct parallel projection method, and it takes $a(k)=\Theta(dk\log{k})$, the runtime on input of size $k$ and $m(k)=k/5$, the average size of the subproblem. Then, the straightforward application of Theorem 2 in \cite{lawlor2013adaptive} gives the following theorem. 

\begin{theorem}
	Assume $N \geq 5k$ and there is no worst-case scenario. Let $T(k)$ denote the runtime of Algorithm \ref{Algorithm2} on a random signal setting with $d_1=1$ and $d_2=d$. Then $\mathbb{E}[T(k)]=\Theta(dk\log{k})$ and 
	\begin{equation*}
	\mathbb{P}[T(k)>\Theta(dk\log{k}) +t dk\log{k}] \leq 5^{-t}.
	\end{equation*}
\end{theorem}

In a similar manner, Algorithm \ref{Algorithm2} takes $a(k)=\Theta(dk)$, the number of samples on input of size $k$ and $m(k)=k/5$, the average size of the subproblem. Thus, Theorem 2 in \cite{lawlor2013adaptive} again gives,

\begin{theorem}
	Assume $N \geq 5k$ and there is no worst-case scenario. Let $S(k)$ denote the number of samples used in  Algorithm \ref{Algorithm2} on a random signal setting with $d_1=1$ and $d_2=d$. Then $\mathbb{E}[S(k)]=\Theta(dk)$ and 
	\begin{equation*}
	\mathbb{P}[S(k)>\Theta(dk)+t dk] \leq 5^{-t}.
	\end{equation*}
\end{theorem}

\section{EMPIRICAL RESULT}\label{5}
\setcounter{equation}{0}

 The partial unwrapping method is implemented in the C language. The pseudo code of this algorithm is shown in Algorithm \ref{Algorithm2}. It is explained in detail in Section \ref{5.0}. In our experiment, dimension $d$ is set to $100$ and $1000$, $d_1$ is $5$ and $d_2$ is $20$ and $200$, accordingly. Frequency bandwidth $N$ in each dimension is $20$ and sparsity $k$ varies as $1, 2, 2^2, \cdots, 2^{10}$. The value of $\epsilon$ for shifting is set to $1/2N^{d_1}$ and the constant number $c$ determining the prime number $p$ is set to $5$.

 We randomly choose $k$ frequency vectors ${\mathbf w}_j \in  \Big[-\frac{N}{2},\frac{N}{2}\Big )^d\cap{\mathbb Z}^d$ and corresponding coefficients $\ a_j=e^{2\pi i\theta_j}\in {\mathbb C}$ from randomly chosen angles $\theta_j\in [0, 1)$ so that the magnitude of each $a_j$ is $1$. For each $d$ and $k$ we have $100$ trials.  We get the result by averaging $\ell^2$ errors, the number of samples used and CPU TICKS out of $100$ trials.

\begin{figure}[ht]
	\centering
  \includegraphics[width=0.90\textwidth, angle=0]{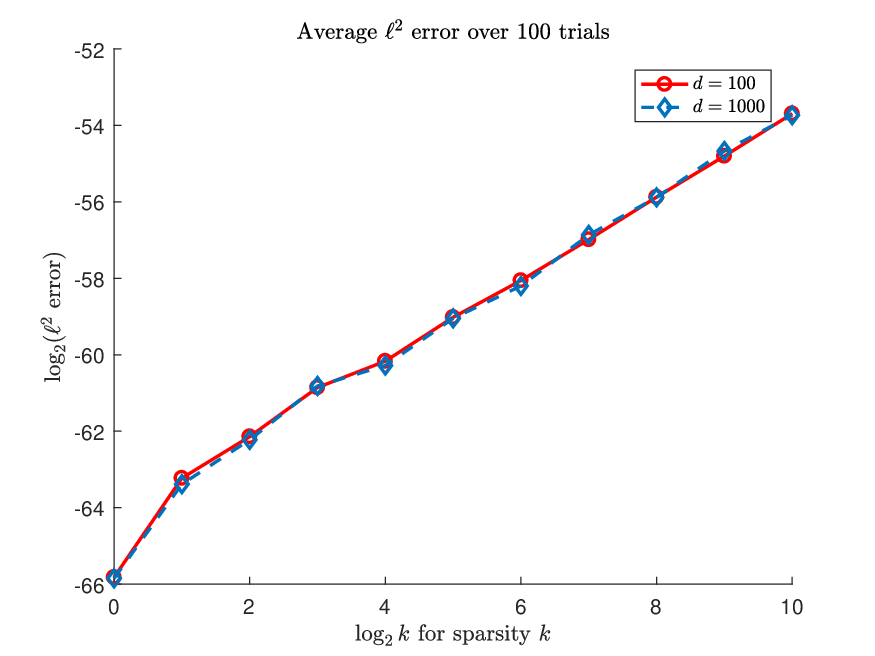}
	\caption{Average $\ell^2$ error}
	\label{fig4}
\end{figure}

\begin{figure}[ht]
	\centering
  \includegraphics[width=0.80\textwidth, angle=0]{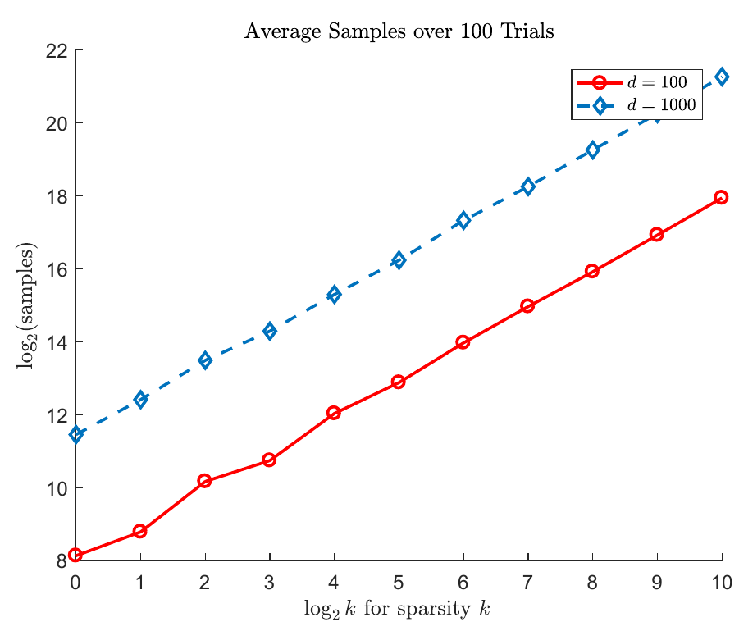}
	\caption{Average sampling complexity}
	\label{fig5}
\end{figure}

\begin{figure}[ht]
	\centering
  \includegraphics[width=0.90\textwidth, angle=0]{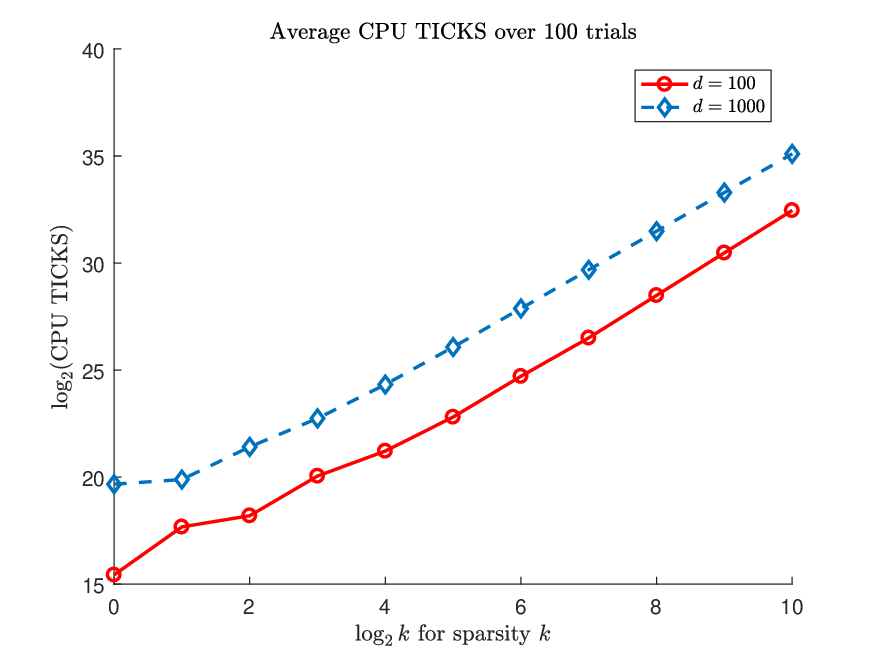}
	\caption{Average CPU TICKS}
	\label{fig6}
\end{figure}

 Since it is difficult to implement high dimensional FFT and there is no practical high dimensional sparse Fourier transform with wide range of $d$ and $s$ at the same time it is hard to compare the result of ours with others, as so far no one else was able to do FFT on this large data set. Thus we cannot help but show ours only. From Figure \ref{fig4} we can see that the average $\ell^2$ errors are below $2^{-52}$. Those errors are from all differences of frequency vectors and coefficients of the original and recovered values. Since all frequency components are integers and thus the least difference is $1$, we can conclude that our algorithm recover the frequency vectors perfectly. Those errors are only from the coefficients. In Figure \ref{fig5} the average sampling complexity is shown. We can see that the logarithm of the number of samples is almost proportional to that of sparsity. Note that the traditional FFT would show the same sampling complexity even though sparsity $k$ varies since it only depends on the bandwidth $N$ and dimension $d$. In Figure \ref{fig6} the average CPU TICKS  are shown. We can see the the logarithm of CPU TICKS is also almost proportional to that of sparsity. Note that the traditional FFT might show the same CPU TICKS even though sparsity $k$ varies since it also depends on the bandwidth $N$ and dimension $d$ only.

\subsection{Algorithm}\label{5.0}
 In this section, the explanation of Algorithm \ref{Algorithm2} is given. In \cite{lawlor2013adaptive} several versions of 1D algorithms are shown. Among them, non-adaptive and adaptive algorithms are introduced where the input function $f$ is not modified throughout the  whole iteration, and is modified by subtracting the function constructed from the data in registry $R$, respectively. In our multidimensional algorithm, however, the adaptive version is mandatory since excluding the contribution of the currently recovered data is the key of our algorithm to avoid the collision of frequencies through projections, whose simple pictorial description is given in Figure \ref{2Dbasic}. In Algorithm \ref{Algorithm2}, the function $g$ is the one constructed from the data in the registry $R$.

 Our algorithm begins with entering inputs, a function $f$, a constant number $c$ determining $p$, a sparsity $k$, a bandwidth $N$ of each dimension, a dimension $d$, factors $d_1$ and $d_2$ of $d$ and a shifting number $\epsilon<1/N$. For each iteration of the algorithm, the number of frequencies to find is updated as $k^{\ast} = k-|R|$. It stops when $|R|$ becomes equal to the sparsity $k$. The prime number $p$ is determined depending on this new $k^{\ast}$ as $p\geq ck^{\ast}$ and is chosen as the next larger prime number. The lines 13 and 14 of Algorithm \ref{Algorithm2} represent the partial unwrapping and sampling with and without shifting from the function where the contribution of former data is excluded. After applying the FFT on each sequence, sorting them according to the magnitude of $\mathcal{F}({f}^{m,n}_{p, 0})$, we  check the ratio between the FFT's of the unshifted and shifted sequences to determined whether there is a collision, either from modulo $p$ or a parallel projection. If all tests are passed, then we find each frequency component and corresponding coefficient for the data that passed and store them in $R$. After several iterations, we find all the data and the final wrapping process gives the original frequency vectors in $d$ dimensions.

\subsection{Accuracy}\label{5.1}
 We assume that there is no noise on the data that we want to recover. Figure \ref{fig4} shows that we can find frequencies perfectly and the $\ell^2$ error from coefficients are significantly small. This error is what we average out over 100 trials for each $d=100, 1000$ and $k=1,2^1,2^2,\cdots,2^{10}$ when $N$ is fixed to $20$. The horizontal axis represents the logarithm with base $2$ of $k$ and the vertical axis represents the logarithm with base $2$ of the $\ell^2$ error. It is increasing as the sparsity $k$ is increasing since the number of nonzero coefficients increases. The red graph in the Figure \ref{fig4} shows the error when the number of dimensions is $100$ and the blue one shows the error when the number of dimensions is $1000$. Thus, we see that the errors are not substantially impacted  by the dimensions.

\subsection{Sampling Complexity}\label{5.2}
 Figure \ref{fig5} shows the sampling complexity of our algorithm averaged out from 100 tests for each dimension and sparsity. The horizontal axis means the logarithm with base $2$ of $k$ and the vertical axis represents the logarithm with base $2$ of the total number of samples from the randomly constructed function which are used to find all frequencies and coefficients. The red graph in the Figure \ref{fig5} shows the sampling complexity when the number of dimensions is $100$ and the blue one shows the one when the number of dimensions is $1000$. Both graphs increase as $k$ increases.  When $d$ is large, we see that it requires more samples since there are more frequency components to find. From the graphs, we  see that the scaling seems to be proportional to $d$.

\subsection{Runtime Complexity}\label{5.3}
 In Figure \ref{fig6}, we plot the runtime complexity of the main part of our algorithm averaged over 100 tests for each dimension and sparsity. \lq\lq Main part'' means that we have excluded the time for constructing a function consisting of frequencies and coefficients and the time associated with getting samples from it. The horizontal axis is the logarithm, base $2$, of $k$ and the vertical axis is the logarithm, base $2$, of CPU TICKS. The red curve shows the runtime when we set the number of dimensions to 100 and the blue one shows the same thing when the number of dimensions to 1000. Both plots increase as $k$ increases. When $d$ is larger, the plots show that it takes more time to run the algorithm. From the graphs we see that the runtime looks proportional to $d$.

 Unfortunately, the sampling process of getting the samples from continuous functions dominates the runtime of the whole algorithm instead of the main algorithm. To show the runtime of our main algorithm, however, we showed CPU TICKS without sampling process. Reducing the time for sampling is still a problem. In \cite{iwen2010combinatorial} the fully discrete Fourier transform is introduced that we expect to use to reduce it. Exploring how to use this will be one part of our future work.

\section{CONCLUSION} \label{6}

In this paper we show how to extend our deterministic $1D$ sublinear sparse Fourier algorithm to the general $d$ dimensional case. The method projects $d$ dimensional frequency vectors onto lower dimensions. In this process we encounter several obstacles. Thus we introduced \lq\lq tilting method'' for the worst case problems and the \lq\lq partial unwrapping method'' to reduce the chance of collisions and to increase the frequency bandwidth within the limit of computation. In this way we can overcome the obstacles as well as maintain the advantage of the $1D$ algorithm. In \cite{lawlor2013adaptive} the sampling complexity is $\Theta(k)$ and the runtime complexity is $\Theta(k{\log}k)$ on average case. Extended this estimation from our $1D$ algorithm, we achieve $\Theta(dk)$ sampling complexity and  a runtime complexity of $\Theta(dk{\log}k)$ on average under the assumption that the worst case scenario does not happen. We will address the complexity of the worst-case scenario in our future work.

 Multidimensional sparse Fourier algorithms have not been discussed much so far, so there is a lot of room for future work. The algorithms in this paper are for recovering data from a noiseless environment only. However most of the actual data contains noise. Thus, the next step will be developing an algorithm for noisy multidimensional data. As mentioned in the previous section, reducing sampling time is another problem to consider. Furthermore, algorithms for fully discrete or nonuniform data will be explored, which is expected to be possible by exploring the way of combining the result from our paper and the one-dimensional fully discrete sparse Fourier transforms from \cite{merhi2017new}. In the end, it is expected that we apply them to real problems like astrophysical data or MRI data.
 
{\bf ACKNOWLEDGEMENTS}
  We would like to thank Mark Iwen for his valuable advice. This research is supported in part by AFOSR grants FA9550-11-1-0281, FA9550-12-1-0343 and FA9550-12-1-0455, NSF grant DMS-1115709, and MSU Foundation grant SPG-RG100059, as well as Hong Kong Research Grant Council grants 16306415 and 16317416.

%\clearpage
\bibliographystyle{abbrv}
\bibliography{phaseshift_multi_d_paper}

\end{document}